# Représentations p-adiques et normes universelles

I. Le cas cristallin

# Bernadette Perrin-Riou

Résumé : Soit $V$ une représentation $p$-adique cristalline du groupe de Galois absolu $G_K$ d'une extension finie non ramifiée $K$ de $\mathbb{Q}_p$ et $T$ un réseau stable par $G_K$. On démontre le résultat suivant :

Soit $\mathrm{Fil}^1 V$ la plus grande sous-représentation de $V$ ayant ses poids de Hodge-Tate strictement positifs et $\mathrm{Fil}^1 T = T \cap \mathrm{Fil}^1 V$. Alors, la limite projective des $H_g^1(K(\mu_{p^n}), T)$ est égale modulo torsion à la limite projective des $H^1(K(\mu_{p^n}), \mathrm{Fil}^1 T)$. En particulier, son rang sur l'algèbre d'Iwasawa est égal à $[K:\mathbb{Q}_p]\dim\mathrm{Fil}^1 V$.

Abstract : Let $V$ be a crystalline $p$-adic representation of the absolute Galois group $G_K$ of an finite unramified extension $K$ of $\mathbb{Q}_p$ and $T$ a lattice of $V$ stable by $G_K$. We prove the following result :

Let $\mathrm{Fil}^1 V$ be the maximal sub-representation of $V$ with Hodge-Tate weights strictly positive and $\mathrm{Fil}^1 T = T \cap \mathrm{Fil}^1 V$. Then, the projective limit of the $H_g^1(K(\mu_{p^n}), T)$ is equal up to torsion to the projective limit of the $H^1(K(\mu_{p^n}), \mathrm{Fil}^1 T)$. So its rank over the Iwasawa algebra is $[K:\mathbb{Q}_p]\dim\mathrm{Fil}^1 V$.

Mathématiques, Bât 425
Université Paris-Sud,
F-91405 Orsay Cedex
France
mél : bpr@geo.math.u-psud.fr

Soit $p$ un nombre premier impair, $K = \mathbb{Q}_p$ et soient $K_n = \mathbb{Q}_p(\mu_{p^n})$ le corps des racines $p^n$-ièmes de l'unité d'ordre $p^n$ dans une clôture algébrique $\overline{\mathbb{Q}}_p$ de $\mathbb{Q}_p$, $K_\infty$ la réunion des $K_n$ et $G_\infty = \text{Gal}(K_\infty/K)$. Soit une courbe elliptique $E$ définie sur $\mathbb{Q}_p$. Si $E(K_n)$ est le groupe des points de $E$ définis sur $K_n$, la limite projective des $E(K_n)$ pour les applications de trace a d'abord été étudiée lorsque $E$ a bonne réduction ordinaire par Mazur et Manin ([9],[10]), puis lorsque $E$ a bonne réduction supersingulière ([6],[16],[11]). Dans le premier cas, c'est un $\mathbb{Z}_p[[G_\infty]]$-module de rang 1, dans le second cas, il est nul. De manière équivalente, dans le cas supersingulier, par la dualité locale de Tate, le groupe de cohomologie galoisienne $H^1(K_\infty, E(\overline{\mathbb{Q}}_p)) = H^1(\text{Gal}(\overline{\mathbb{Q}}_p/K_\infty), E(\overline{\mathbb{Q}}_p))$ est nul et en notant $E_{p^\infty}$ le groupe des points de torsion $p$-primaire de $E(\overline{\mathbb{Q}}_p)$, on a alors l'isomorphisme

$$\mathbb{Q}_p/\mathbb{Z}_p \otimes E(K_\infty) \cong H^1(K_\infty, E_{p^\infty}) \ .$$

Des résultats de ce type ont été étendus aux groupes formels ([16], [11]). Soit maintenant $V$ une représentation $p$-adique de $G_{\mathbb{Q}_p} = \text{Gal}(\overline{\mathbb{Q}}_p/\mathbb{Q}_p)$, c'est-à-dire un $\mathbb{Q}_p$-espace vectoriel de dimension finie muni d'une action linéaire et continue de $G_{\mathbb{Q}_p}$ et soit $T$ un $\mathbb{Z}_p$-réseau de $V$ stable par $G_{\mathbb{Q}_p}$. L'analogue des $E(K)$ ou plutôt du complété $p$-adique $\varprojlim_n E(K)/p^n E(K)$ a été fourni par Bloch et Kato : pour $* \in \{e, f, g\}$, ils construisent des sous-$\mathbb{Q}_p$-espaces vectoriels $H^1_*(K, V)$ du groupe de cohomologie galoisienne $H^1(K, V) = H^1(\text{Gal}(\overline{\mathbb{Q}}_p/K), V)$ et un $\mathbb{Z}_p$-module de type fini $H^1_*(K, T)$ comme l'image réciproque de $H^1_*(K, V)$ dans $H^1(K, T)$.

La question suivante se pose alors : que peut-on dire de la limite projective des $H^1_*(K_n, T)$ que nous noterons $Z_{\infty,*}(K, T)$ ou, de manière duale, que peut-on dire de la différence entre $H^1(K_\infty, V/T)$ et $H^1_*(K_\infty, V/T)$ où $H^1_*(K_\infty, V/T)$ est la limite inductive des $\mathbb{Q}_p/\mathbb{Z}_p \otimes H^1_*(K_n, T)$ ? C'est ce problème que nous abordons ici. Nous prenons maintenant pour $K$ une extension finie de $\mathbb{Q}_p$.

Nous supposons désormais que $V$ est une représentation $p$-adique de Hodge-Tate, c'est-à-dire qu'il existe des entiers $j$ tels que, si $\mathbb{C}_p$ est le complété $p$-adique de $\overline{\mathbb{Q}}_p$, on a la décomposition de Hodge-Tate de $G_K$-modules

$$\mathbb{C}_p \otimes V \cong \oplus_{i \in \mathbb{Z}} \mathbb{C}_p \otimes_K (\mathbb{C}_p \otimes V(-i))^{G_K}(i)$$

où $W(i)$ désigne le twist à la Tate $W \otimes \mathbb{Q}_p(i)$ (l'action de Galois est twistée par la puissance $i$-ième du caractère cyclotomique $\chi : G_\infty \to \mathbb{Z}_p^*$). Les nombres de Hodge-Tate sont alors les $h_j(V) = \dim_K(\mathbb{C}_p \otimes V(-i))^{G_K}$, les poids de Hodge-Tate sont les entiers $j$ tels que $h_j(V)$ est non nul.

Les trois $\mathbb{Z}_p[[G_\infty]]$-modules $Z_{\infty,*}(K, T)$ sont essentiellement les mêmes. Rappelons que $Z_{\infty,g}(K, T)$ est un sous-$\mathbb{Z}_p[[G_\infty]]$-module de la limite projective $Z^1_\infty(K, T)$ des modules $H^1(K_n, T)$, que le module de torsion de $Z^1_\infty(K, T)$ (en tant que $\mathbb{Z}_p[[\Gamma]]$-module avec $\Gamma = \text{Gal}(K_\infty/K_1)$) est isomorphe à $T^{G_{K_\infty}}$. Le sous-module de torsion de $Z_{\infty,g}(K, T)$ se calcule aisément. Aussi, nous nous intéressons au quotient $\tilde{Z}_{\infty,g}(K, V)$ de $Z_{\infty,g}(K,V) = \mathbb{Q}_p \otimes Z_{\infty,g}(K,T)$ par son $\mathbb{Q}_p \otimes \mathbb{Z}_p[[\Gamma]]$-module de torsion et même plus simplement au rang de $Z_{\infty,g}(K,T)$ comme $\mathbb{Z}_p[[G_\infty]]$-module ou plutôt au $\mathbb{Z}_p[[\Gamma]]$-rang des sous-espaces propres sous l'action de $\Delta = \mathrm{Gal}(K_1/K)$.

Lorsque les poids de Hodge-Tate de $V$ sont tous strictement positifs, on a $H^1_g(K_n, V) = H^1(K_n, V)$. Donc, $Z_{\infty,g}(K,T)$ est égal à $Z^1_\infty(K,T) = \varprojlim_n H^1(K_n, T)$ et est donc un $\mathbb{Z}_p[[G_\infty]]$-module de rang $[K:\mathbb{Q}_p]\dim_{\mathbb{Q}_p} V$. Lorsque les poids de Hodge-Tate sont négatifs ou nuls, l'espace tangent de $V$ est nul et $Z_{\infty,g}(K,T)$ est de torsion. Dans ces deux cas extrêmes, le problème posé est donc vite résolu.

Le premier cas qui a été étudié est celui où la représentation $V$ est de Dabrowski-Panschiskhin (D-P), c'est-à-dire qu'il existe une sous-représentation $p$-adique $\mathrm{Fil}^1 V$ de $V$ telle que les poids de Hodge-Tate de $\mathrm{Fil}^1 V$ sont tous strictement positifs et les poids de Hodge-Tate de $V/\mathrm{Fil}^1 V$ soient tous négatifs on nuls. C'est en particulier le cas lorsque $V$ est une représentation $p$-adique ordinaire, avec le cas encore plus particulier où $V$ a un seul poids de Hodge-Tate.

**0.1. Proposition.** *Soit $V$ une représentation $p$-adique de D-P telle que $V^{G_{K_\infty}} = 0$. Alors, $Z_{\infty,g}(K,V)$ est isomorphe à $Z^1_\infty(K, \mathrm{Fil}^1 V)$. C'est un $\mathbb{Q}_p \otimes \mathbb{Z}_p[[G_\infty]]$-module de rang $[K:\mathbb{Q}_p]\dim_{\mathbb{Q}_p} \mathrm{Fil}^1 V = [K:\mathbb{Q}_p] \sum_{j>0} h_j(V)$.*

Remarquons que $[K:\mathbb{Q}_p] \sum_{j>0} h_j(V)$ est aussi la dimension de l'espace tangent $t_V(K)$ de $V$ sur $K$.

Lors d'une discussion où je lui parlais d'un résultat partiel dans le cas où $V$ est irréductible, Nekovář a proposé la conjecture suivante :

**Conjecture A.** *Soit $V$ une représentation $p$-adique de $G_K$ qui est de de Rham. Soit $\mathrm{Fil}^1 V$ la plus grande sous-représentation $p$-adique de $V$ dont les poids de Hodge-Tate sont strictement positifs. Alors, $Z^1_\infty(K, \mathrm{Fil}^1 V)$ est contenu dans $Z_{\infty,g}(K,V)$ et le conoyau est un $\mathbb{Z}_p[[\Gamma]]$-module de torsion.*

**Définition.** On dit que $V$ vérifie (N) si $V$ n'admet pas de sous-représentation non triviale $W$ telle que les nombres de Hodge-Tate de $W$ soient tous strictement positifs.

Le quotient $V/\mathrm{Fil}^1 V$ vérifie la propriété (N). Par un argument de dévissage, on montre facilement que la conjecture (A) est équivalente à la conjecture suivante :

**Conjecture A'.** *Soit $V$ une représentation $p$-adique de de Rham vérifiant (N). Alors $Z_{\infty,g}(K,V)$ est un $\mathbb{Q}_p \otimes \mathbb{Z}_p[[\Gamma]]$-module de torsion.*

Nous démontrons ici cette conjecture lorsque $V$ est une représentation $p$-adique cristalline.

**0.2. Théorème.** *Soit $K$ une extension finie non ramifiée de $\mathbb{Q}_p$. Si $V$ est une représentation $p$-adique cristalline de $G_K$, $Z^1_\infty(K, \mathrm{Fil}^1 V)$ est contenu dans $Z_{\infty,g}(K,V)$ et le quotient est un $\mathbb{Z}_p[[\Gamma]]$-module de torsion contenu dans $(V/\mathrm{Fil}^1 V)^{G_{K_\infty}}/(V/\mathrm{Fil}^1 V)^{G_K}$.*

On a ainsi une suite exacte
$$0 \to Z^1_\infty(K, \mathrm{Fil}^1 V) \to Z_{\infty,g}(K,V) \to (V/\mathrm{Fil}^1 V)^{G_{K_\infty}}/(V/\mathrm{Fil}^1 V)^{G_K}.$$

Remarquons que les $H^1_g(K_n, V)$ sont à peu près de dimension égale à la dimension de $t_V(K_n) = [K_n : K]t_V(K)$. Par contre, $\text{Fil}^1 V$ est en général de dimension bien inférieure à la dimension de $t_V(K)$ et le rang de $Z_{\infty,g}(K, V)$ est donc en général loin d'être maximal. Le cas D-P est en fait très particulier.

Nous démontrerons ce théorème sous la forme équivalente :

**0.3. Théorème.** *Soit $K$ une extension finie non ramifiée de $\mathbb{Q}_p$. Si $V$ est une représentation $p$-adique cristalline de $G_K$ vérifiant (N), alors $Z_{\infty,g}(K, V)$ est de torsion sur $\mathbb{Q}_p \otimes \mathbb{Z}_p[[\Gamma]]$.*

Autrement dit, si $V$ une représentation $p$-adique cristalline vérifiant (N) et telle que $V^{G_{K_\infty}} = 0$ :

1. Il n'existe pas de famille compatible (pour les applications trace) d'éléments de $H^1_g(K_n, T)$ ;
2. par dualité,
$$H^1(K_\infty, V^*(1)/T^*(1)) = H^1_g(K_\infty, V^*(1)/T^*(1))$$

**0.4. Corollaire.** *Si $V$ est une représentation $p$-adique cristalline telle que $V^{G_{K_\infty}} = 0$, le sous-$\mathbb{Z}_p$-module de $H^1_g(K, T)$ des normes (traces) universelles d'éléments de $H^1_g(K_n, T)$ relativement à l'extension $K_\infty/K$ est de rang $[K : \mathbb{Q}_p] \dim \text{Fil}^1 V$.*

Par dualité l'énoncé 0.2 devient :

**0.5. Corollaire.** *Soit $K$ une extension finie non ramifiée de $\mathbb{Q}_p$ et soit $V$ une représentation $p$-adique cristalline de $G_K$ telle que $V^{G_{K_\infty}} = 0$.*

*1) Si $V$ n'admet pas de représentation-quotient dont les poids de Hodge-Tate soient tous négatifs, alors $H^1_g(K_\infty, V/T) = \varinjlim H^1_g(K_n, V/T) = 0$.*

*2) Soit $W$ le plus grand quotient de $V$ dont les poids de Hodge-Tate sont tous négatifs et $U$ l'image de $T$ dans $W$. Alors,*
$$H^1_g(K_\infty, V/T) = H^1_g(K_\infty, W/U) .$$

En particulier, si $V$ est la représentation $p$-adique associée à une courbe elliptique, on retrouve que $Z_{\infty,g}(\mathbb{Q}_p, V)$ est nul lorsque $V$ est cristalline (c'est-à-dire que $E$ a bonne réduction en $p$) et irréductible (c'est-à-dire que $E$ a réduction supersingulière). Plus généralement, soit $V_f$ la restriction à $\mathbb{Q}_p$ de la représentation $p$-adique sur $\mathbb{Q}$ attachée à une forme modulaire $f = \sum_{n>0} a_n(f)q^n$ de poids $k$ de niveau premier à $p$. C'est une représentation $p$-adique cristalline de dimension 2. Pour préciser $V_f$, disons que ses poids de Hodge-Tate sont 0 et $k-1$. La représentation $V_f$ est irréductible si et seulement si elle n'est pas ordinaire, c'est-à-dire si et seulement si $a_p(f)$ n'est pas une unité. On déduit du théorème 0.2 le théorème suivant :

**0.6. Théorème.** *1) Si $a_p(f)$ est divisible par $p$, $Z_{\infty,g}(K, V_f(j))$ est nul pour $j \leq 0$ et de $\mathbb{Z}_p[[G_\infty]]$-rang 2 pour $j > 0$.*

*2) Si $a_p(f)$ n'est pas divisible par $p$, $Z_{\infty,g}(K, V_f(j))$ est un $\mathbb{Z}_p[[G_\infty]]$-module de rang 0 pour $j \leq -k+1$, de rang 1 pour $-k+2 \leq j \leq 0$ et de rang 2 pour $j > 0$.*

En utilisant les méthodes standard, on peut déduire de ces résultats des conséquences sur le rang du groupe de Selmer d'une représentation $p$-adique d'un corps de nombres ayant bonne réduction absolue aux places divisant $p$. Nous ne le ferons pas ici.

Le théorème 0.2 implique le théorème suivant

**0.7. Théorème.** *Soient $K$ une extension finie non ramifiée de $\mathbb{Q}_p$ et $V$ une représentation cristalline $p$-adique de $G_K$. Par l'isomorphisme naturel de twist*

$$Z^1_\infty(K,V) \to Z^1_\infty(K,V(1)) ,$$

*l'image de $Z_{\infty,g}(K,V)$ est contenue dans $Z_{\infty,g}(K,V(1))$.*

On peut bien sûr conjecturer que cela reste vrai lorsque $V$ est une représentation de de Rham.

**Remarques.** Une conséquence facile par dévissage du théorème 0.2 est que la conjecture (A) est vraie pour une représentation $p$-adique de de Rham, extension d'une extension cristalline par une représentation de de Rham de poids de Hodge-Tate strictement positifs. Par exemple, elle est vraie pour une extension semi-stable $V$ de $\mathbb{Q}_p$ par $\mathbb{Q}_p(1)$ (ce que l'on savait déjà puisqu'il s'agit d'une représentation ordinaire). On peut cependant construire des exemples un peu moins triviaux :

Soit $W$ une représentation cristalline. Il existe une extension $V$ de $ad(W)^0$ par $\mathbb{Q}_p(1)$ qui ne soit pas cristalline. En effet les extensions de $ad(W)^0$ par $\mathbb{Q}_p(1)$ sont classifiées par $H^1(K, ad(W)^0 \otimes \mathbb{Q}_p(1))$, les extensions cristallines par $H^1_f(K, ad(W)^0 \otimes \mathbb{Q}_p(1))$ et les extensions semi-stables par $H^1_g(K, ad(W)^0 \otimes \mathbb{Q}_p(1))$. Comme $p^{-1}$ est valeur propre de $\varphi$ sur $\mathbf{D}_p(ad(W)^0 \otimes \mathbb{Q}_p(1))$, il existe une telle extension $V$ qui soit semistable et non cristalline. Le théorème est alors vrai pour $V$. Calculons $\mathrm{Fil}^1 V$. Supposons que $ad(W)^0$ est irréductible. Il est clair que $\mathrm{Fil}^1 V$ contient $\mathbb{Q}_p(1)$. S'il ne lui est pas égal, $\mathrm{Fil}^1 V/\mathbb{Q}_p(1)$ a des poids de Hodge-Tate $> 0$ et l'extension $\mathrm{Fil}^1 V$ de $\mathrm{Fil}^1 V/\mathbb{Q}_p(1)$ par $\mathbb{Q}_p(1)$ est nécessairement triviale. Ce qui contredit l'irréductibilité de $ad(W)^0$. Ainsi, $\mathrm{Fil}^1 V = \mathbb{Q}_p(1)$ et $Z_{\infty,g}(K,V)$ est de rang $[K:\mathbb{Q}_p]$ et égal à $Z_{\infty,g}(K,\mathbb{Q}_p(1))$ modulo torsion. D'autre part, la représentation $V$ n'est pas de D-P car $V/\mathrm{Fil}^1 V$ a des poids positifs et négatifs.

Disons un mot sur la démonstration. En utilisant la théorie de [13], on obtient une description complète de $Z_{\infty,g}(K,V)$ en termes de fonctions analytiques. On est alors ramené à un problème de nullité de fonctions analytiques vérifiant certaines propriétés. Pour cela, on majore leur ordre et on minore le "nombre de zéros" afin d'aboutir à une contradiction.

Dans la première partie, nous rappelons quelques faits simples et bien connus sur $Z^1_\infty(K,V)$ et nous démontrons les implications reliant les théorèmes énoncés dans l'introduction. Dans la deuxième partie nous introduisons les techniques nécessaires à la démonstration du théorème 0.3. Dans la troisième partie, nous traitons le cas où $V$ est de dimension 2. Dans la quatrième partie, nous traitons le cas général.

Nous prendrons maintenant pour $K$ une extension finie non ramifiée de $\mathbb{Q}_p$ et $K_n = K(\mu_{p^n})$. On note $\sigma$ l'homomorphisme de Frobenius de $K$. On fixe une représentation $p$-adique $V$ de $G_K$ qui soit cristalline.

## 1. Théorie d'Iwasawa classique

### 1.1. Sous-modules de torsion.
L'application d'inflation

$$(T^{G_{K_\infty}})_{\mathrm{Gal}(K_\infty/K_n)} \cong H^1(K_\infty/K_n, T^{G_{K_\infty}}) \to H^1(K_n, T)$$

induit un homomorphisme injectif de $T^{G_{K_\infty}}$ dans $Z^1_\infty(K,T)$. On démontre que $T^{G_{K_\infty}}$ est le module de $\mathbb{Z}_p[[\Gamma]]$-torsion de $Z^1_\infty(K,T)$ ([12]).

Comme $V$ est de Hodge-Tate et $K$ non ramifiée sur $\mathbb{Q}_p$, $V^{G_{K_\infty}}$ est isomorphe à $\oplus_i (V^{G_{K_\infty}}(-i))^{G_K}(i)$ ([13, 3.4.3]).

**1.1.1. Proposition.** *Le $Zp[[\Gamma]]$-module de torsion de $Z_{\infty,g}(K,T)$ est isomorphe à $V^{G_{K_\infty}}/V^{G_K} = \oplus_{i\neq 0}(V^{G_{K_\infty}}(-i))^{G_K}(i)$.*

*Démonstration.* Par définition de $H^1_g(K_n,T)$ comme image réciproque de $H^1_g(K_n,V)$ dans $H^1(K_n,T)$, $H^1_g(K_n,T)$ contient le sous-module de $\mathbb{Z}_p$-torsion de $H^1(K_n,T)$. Comme $H^1(K_\infty/K_n, \oplus_{i\neq 0}(V^{G_{K_\infty}}(-i))^{G_K}(i)) = 0$, l'image de $\oplus_{i\neq 0}(T^{G_{K_\infty}}(-i))^{G_K}(i))$ dans $Z^1_\infty(K,T)$ est contenue dans $Z_{\infty,g}(K,T)$. Il reste à regarder l'image de $T^{G_K} \cong \varprojlim H^1(K_\infty/K_n, T^{G_K})$. On est alors ramené au cas de la représentation triviale et à démontrer que l'application $H^1(K_\infty/K_n, \mathbb{Q}_p) \to H^1(K_n, \mathbb{Q}_p)/H^1_g(K_n, \mathbb{Q}_p)$ est injective. Or $H^1_g(K_n, \mathbb{Q}_p)$ s'interprète comme $\mathrm{Hom}_{\mathbb{Z}_p}(\mathrm{Gal}(K^{nr}/K), \mathbb{Q}_p)$ (où $K^{nr}$ est la plus grande extension non ramifiée de $K$) qui est d'intersection nulle avec $\mathrm{Hom}_{\mathbb{Z}_p}(\mathrm{Gal}(K_\infty/K_n), \mathbb{Q}_p)$ car $K_\infty$ est totalement ramifiée. □

1.2. Montrons que le théorème 0.3 implique le théorème 0.2. Soit $V' = \mathrm{Fil}^1 V$ comme dans l'introduction et $V'' = V/V'$. La représentation $V''$ vérifie la la condition (N). Donc $Z_{\infty,g}(K,V'')$ est de torsion. D'autre part, les poids de Hodge-Tate de $V'$ sont strictement positifs. On en déduit que $Z_{\infty,g}(K,V') = Z_\infty(K,V')$ (on a a en effet $H^1_g(K_n,T') = H^1(K_n,T')$ pour tout entier $n$). D'autre part, comme $V$ est de de Rham, la suite exacte $0 \to T' \to T \to T'' \to 0$ induit la suite exacte

$$0 \to\to T'^{G_{K_n}} \to T^{G_{K_n}} \to T'^{G_{K_n}} \to H^1_g(K_n,T') \to H^1_g(K_n,T) \to H^1_g(K_n,T'')$$

Par passage à la limite projective, les trois premiers termes sont nuls (pour $n$ assez grand, les modules du type $T^{G_{K_n}}$ sont des $\mathbb{Z}_p$-modules de type fini stationnaires et l'application de corestriction est la multiplication par $p$). On obtient donc la suite exacte

$$0 \to Z_{\infty,g}(K,V') \to Z_{\infty,g}(K,V) \to Z_{\infty,g}(K,V'') .$$

On en déduit que l'image de $Z_\infty(K,V')$ dans $Z_\infty(K,V)$ est contenue dans $Z_{\infty,g}(K,V)$ et que l'on a une suite exacte de $\mathbb{Z}_p[[G_\infty]]$-modules

$$0 \to Z_\infty(K,V') \to Z_{\infty,g}(K,V) \to V''^{G_{K_\infty}}/V''^{G_K} .$$

Le dernier module contient l'image de $V^{G_{K_\infty}}$. Mais la flèche
$$V^{G_{K_\infty}} \to V''^{G_{K_\infty}}/V''^{G_K}$$
peut être ou ne pas être surjective. Je ne sais pas si $\tilde{Z}_\infty(K, \mathrm{Fil}^1 V) = \tilde{Z}_{\infty,g}(K, V)$.

1.3. Montrons que le théorème 0.2 implique le théorème 0.7. On utilise dans ce qui suit les modules filtrés de Fontaine $\mathbf{D}_{\mathrm{dR}}(V)$. Remarquons que $\mathrm{Fil}^0 \mathbf{D}_{\mathrm{dR}}(\mathrm{Fil}^1 V) = 0$. On a alors $(\mathrm{Fil}^1 V)(1) \subset \mathrm{Fil}^1(V(1))$ : en effet, $\mathrm{Fil}^0 \mathbf{D}_{\mathrm{dR}}(\mathrm{Fil}^1 V(1)) = \mathrm{Fil}^1 \mathbf{D}_{\mathrm{dR}}(\mathrm{Fil}^1 V) \subset \mathrm{Fil}^0 \mathbf{D}_{\mathrm{dR}}(\mathrm{Fil}^1 V) = 0$. L'image de $Z_{\infty,g}(K, V) = Z_\infty^1(K, \mathrm{Fil}^1 V)$ est $Z_\infty^1(K, \mathrm{Fil}^1 V(1))$ qui est contenue dans $Z_\infty^1(K, \mathrm{Fil}^1(V(1))) = Z_{\infty,g}(K, V(1))$ (à torsion près).

On peut conjecturer que la conclusion du théorème 0.7 reste vraie si l'on suppose seulement que $V$ est une représentation $p$-adique de de Rham. On montre de même que c'est une conséquence de la conjecture A. Réciproquement, si le théorème 0.7 est vrai pour les représentations de de Rham et si la conjecture A est vraie pour les représentations de de Rham de poids de Hodge-Tate positifs ou nuls, la conjecture A est vraie pour toute représentation $p$-adique de de Rham : On peut supposer que $\mathrm{Fil}^1 V = 0$. Soit $u^* - 1$ le plus grand poids de Hodge-Tate de $\mathbf{D}_{\mathrm{dR}}(V)$. On peut supposer que $u^* \geq 1$. Si $u^* = 1$, la conjecture A est supposée vraie. Montrons l'assertion par récurence sur $u^*$. Posons $W = \mathrm{Fil}^1(V(1))(-1)$, c'est une sous-représentation de $V$. L'image de $Z_{\infty,g}(K, V)$ est contenue dans $Z_{\infty,g}(K, V(1)) = Z_\infty^1(K, W(1))$, ce qui signifie que $Z_{\infty,g}(K, V)$ est en fait contenu dans $Z_\infty^1(K, W)$ et même dans $Z_{\infty,g}(K, W)$. Or, $\mathrm{Fil}^1 \mathbf{D}_{\mathrm{dR}}(W) = \mathrm{Fil}^0 \mathbf{D}_{\mathrm{dR}}(\mathrm{Fil}^1(V(1))) = 0$, les poids de Hodge-Tate de $W$ sont donc positifs ou nuls. Comme $\mathrm{Fil}^1 W \subset \mathrm{Fil}^1 V$, $\mathrm{Fil}^1 W$ est nul. On déduit de nouveau que $Z_{\infty,g}(K, W)$ est de torsion. Ce qui démontre la conjecture A.

1.4. **Normes universelles.** Supposons maintenant pour simplifier que $V^{G_{K_\infty}}$ et $(V/\mathrm{Fil}^1 V)^{G_{K_\infty}}$ sont nuls.

L'application naturelle $Z_\infty^1(K, V)_{G_\infty} \to H^1(K, V)$ est une injection ([13, 3.2]). On a alors le diagramme commutatif :

$$\begin{array}{ccc} Z_\infty^1(K, \mathrm{Fil}^1 V)_{G_\infty} = Z_{\infty,g}^1(K, V)_{G_\infty} & \to & Z_\infty^1(K, V)_{G_\infty} \\ \downarrow & & \downarrow \\ H^1(K, \mathrm{Fil}^1 V) & \to & H^1(K, V) \end{array}$$

L'hypthèse $(V/\mathrm{Fil}^1 V)^{G_K} = 0$ implique que la flèche d'en bas est injective, l'hypothèse plus forte que $(V/\mathrm{Fil}^1 V)^{G_{K_\infty}} = 0$ que la flèche d'en haut l'est. On en déduit que le sous-espace de $H^1(K, V)$ formé des normes universelles relativement aux $H_g^1(K_n, T)$ (c'est-à-dire l'image de $Z_{\infty,g}(K, V)$) est un $\mathbb{Q}_p$-espace vectoriel de dimension $[K : \mathbb{Q}_p] \dim \mathrm{Fil}^1 V$.
.

## 2. Préliminaires

2.1. **Un lemme sur les fonctions analytiques.** Soit $\mathcal{H}$ l'algèbre des séries formelles en une variable $x$ à coefficients dans $\mathbb{Q}_p$ convergeant sur le disque unité $\{x \in$ $\mathbb{C}_p$ tel que $|x| < 1\}$ où $\mathbb{C}_p$ est le complété $p$-adique de $\overline{\mathbb{Q}_p}$. Si $\rho$ est un réel inférieur à 1, on note $||f||_\rho = \sup_{|x|=\rho} |f(x)| = \sup_{|x|\leq\rho} |f(x)|$. On a alors ([3, IV])

$$||fg||_\rho = ||f||_\rho ||g||_\rho .$$

Il existe un opérateur continu linéaire $\varphi$ sur $\mathcal{H}$ tel que $\varphi(1+x) = (1+x)^p$. On vérifie facilement que

$$||\varphi(f)||_\rho = ||f||_{\rho^p}$$

pour $\rho > p^{-1/(p-1)}$. En particulier, si $\rho_n = p^{-\frac{1}{p^n(p-1)}}$, on a

$$||\varphi^n(f)||_{\rho_n} = ||f||_{\rho_1} .$$

Si $r \in \mathbb{R}$, on note $\mathcal{H}_r$ le sous-$K$-espace vectoriel de $\mathcal{H}$ formé des séries $F$ telles que la suite $||p^{nr}F||_{\rho_n}$ est bornée. On définit alors $\mathfrak{o}(F)$ (resp. $\mathfrak{O}(F)$) comme la borne inférieure (resp. le plus petit) des réels $r$ tels que la suite $||p^{nr}F||_{\rho_n}$ tende vers 0 (resp. soit bornée) avec $n$. Par exemple, $\mathfrak{O}(\log^r) = r$.

On pose $\mathcal{H}_K = K \otimes \mathcal{H}$ (resp. $\mathcal{H}_{K,r} = K \otimes \mathcal{H}_r$) et on étend $\varphi$ à $\mathcal{H}_K$ par $\sigma$-linéarité.

**2.1.1. Lemme.** *Soient $f$ et $g$ deux éléments de $\mathcal{H}_K$. On suppose que*

$$\frac{\varphi(f)}{f} = \mu \frac{\varphi(g)}{g}$$

*avec $\mu \in K$. Alors, $\mathfrak{O}(f)$, resp. $\mathfrak{o}(f)$ existe si et seulement si $\mathfrak{O}(g)$, resp. $\mathfrak{o}(g)$ existe et on a alors*

$$\mathfrak{O}(f) = \operatorname{ord} \mu + \mathfrak{O}(g) \quad \text{resp. } \mathfrak{o}(f) = \operatorname{ord} \mu + \mathfrak{o}(g) .$$

*Démonstration.* On a

$$\frac{f}{\varphi^n(f)} = (\prod_{i=0}^{n-1} \sigma^i(\mu))^{-1} \frac{f}{\varphi^n(f)} .$$

En prenant la norme $||.||_{\rho_n}$, on en déduit que pour tout entier $n$, on a

$$||f||_{\rho_n} = C p^{n \operatorname{ord} \mu} ||g||_{\rho_n}$$

D'où,

$$||p^{ns} f||_{\rho_n} = C ||p^{n(s-\operatorname{ord} \mu)} g||_{\rho_n}$$

Par définition de $\mathfrak{O}$ et de $\mathfrak{o}$, le lemme s'en déduit. □

Soit $\psi$ l'opérateur de $\mathbb{Q}_p[[x]]$ tel que $\varphi \circ \psi(g) = p^{-1} \sum_{\zeta \in \mu_p} f(\zeta(1+x) - 1)$ prolongé par $\sigma^{-1}$-linéarité à $\mathcal{H}_K$. Le noyau de $\psi$ sur $\mathcal{H}_r$ est canoniquement isomorphe à $\mathcal{H}_r(G_\infty)$ où $\mathcal{H}_r(G_\infty) = \mathbb{Z}_p[\operatorname{Gal}(K_1/K)] \otimes \mathcal{H}_r(\operatorname{Gal}(K_\infty/K_1))$ avec $\mathcal{H}_r(\operatorname{Gal}(K_\infty/K_1)) = \{\sum_n a_n (\gamma-1)^n$ avec $\sum_n a_n X^n \in \mathcal{H}_r\}$ et $\gamma$ un générateur topologique de $\operatorname{Gal}(K_\infty/K_1)$), l'isomorphisme $\mathcal{H}_r(G_\infty) \to \mathcal{H}_r^{\psi=0}$ est induit par $\tau \mapsto \tau.(1+x) = (1+x)^{\chi(\tau)}$, nous le noterons $h \mapsto h.(1+x)$ (pour prolonger à $\mathcal{H}_r(G_\infty)$, on montre que si $f_{n,r}$ est le polynôme d'approximation de $f$ modulo $\prod_{i=0}^r (\chi^{-i}(\gamma_n)\gamma_n - 1)$, la suite $f_{n,r}.(1+x)$ converge dans $\mathcal{H}_r$ et ne dépend pas des choix ; c'est par définition $f.(1+x)$).

L'opérateur de dérivation $D$ sur $\mathcal{H}$ donné par $D(g) = (1+x)g'$ restreint à $\mathcal{H}^{\psi=0}$ correspond sur $\mathcal{H}(G_\infty)$ à l'opération de twist $\tau \mapsto \chi(\tau)\tau$.

Rappelons le lemme suivant ([13, 1.3]) :

**2.1.2. Lemme.** *Soit $f \in \mathcal{H}_r$. Supposons que $f(\zeta - 1) = 0$ pour toute racine de l'unité d'ordre une puissance de $p$. Alors, il existe $g \in \mathcal{H}_{r-1}$ tel que $f = g\log(1+x)$.*

Nous dirons dans ce cas que $f$ est divisible par $\log(1+x)$ (dans $\mathcal{H}$).

## 2.2. Un lemme de déterminant.
Soit $W$ un $K$-espace vectoriel de dimension finie $d$.

**Lemme.** *Soient $g_1, \ldots, g_d$ des éléments de $\mathcal{H}_K \otimes_K W$. On suppose que pour tout entier $n$, il existe une filtration décroissante exhaustive et séparée $\mathrm{Fil}^i W_n$ de $W_n = K_n \otimes_K W$ avec $\mathrm{Fil}^1 W_n = 0$ telle que*

1. *les entiers $h_j = \dim_{K_n} \mathrm{Fil}^j W_n - \dim_{K_n} \mathrm{Fil}^{j+1} W_n$ ne dépendent pas de $n$ ;*
2. *pour tout entier $j \leq 0$ et toute racine de l'unité $\zeta_n$ d'ordre $p^n$,*

$$D^{-j}(g_i)(\zeta_n - 1) \in \mathrm{Fil}^j W_n\ .$$

*Alors $\det(g_1, \cdots, g_d)$ (calculé dans une base de $W$) est divisible par $\log^{-t_H}(1+x)$ où $t_H = \sum_j j h_j$ est le degré de la filtration.*

*Démonstration.* Il est clair que $r_j = \dim \mathrm{Fil}^j W_n$ ne dépend pas non plus de $n$ et on a $h_j = r_j - r_{j+1}$. Posons $F = \det(g_1, \cdots, g_d)$. Fixons un entier $n$ et $\zeta$ une racine de l'unité d'ordre $p^n$. Il s'agit de démontrer que $D^i(F)(\zeta - 1) = 0$ pour tout entier $i$ avec $0 \leq i < r$. Or $D^j(F)(\zeta_n - 1) \in \sum_{\sum_i k_i = j} \wedge_{i=1}^d \mathrm{Fil}^{k_i} W_n = \mathrm{Fil}^j \wedge^d W_n$. Comme $\wedge^d W_n$ a un unique poids de Hodge $t_H$, $D^j(F)$ est divisible par $\log(1+x)$ tant que $j < -t_H$, ce qui implique que $F$ est divisible par $\log^{-t_H}(1+x)$. □

## 2.3. $\varphi$-modules filtrés.
A la représentation $p$-adique cristalline $V$ est associé un $\varphi$-module filtré $\mathbf{D}_p(V)$ sur $K$, c'est-à-dire un $K$-espace vectoriel de dimension $d$, muni d'un opérateur $\sigma$-linéaire bijectif $\varphi$ et d'une filtration $\mathrm{Fil}^i \mathbf{D}_p(V)$ décroissante, exhaustive et séparée. Posons $h_j = h_j(\mathbf{D}_p(V)) = \dim_K \mathrm{Fil}^j \mathbf{D}_p(V) - \dim_K \mathrm{Fil}^{j+1} \mathbf{D}_p(V)$, on a alors $h_j(\mathbf{D}_p(V)) = h_{-j}(V)$. Les poids de Hodge de $\mathbf{D}_p(V)$ sont les entiers $j$ tel que $h_j(\mathbf{D}_p(V)) \neq 0$. Soit $t_H(\mathbf{D}_p(V)) = \sum_j j h_j(\mathbf{D}_p(V))$ le degré de Hodge du $\varphi$-module filtré $\mathbf{D}_p(V)$. On note $-u$ (resp. $u^* - 1$) le plus petit (resp. le plus grand) poids de Hodge de $\mathbf{D}_p(V)$ et on pose $I_H(\mathbf{D}_p(V)) = \{-u, u^* - 1\}$. Si $\Delta$ est un sous-$K$-espace vectoriel de $\mathbf{D}_p(V)$ stable par $\varphi$ (sous-$\varphi$-module), on le munit de la filtration induite (sous-$\varphi$-module filtré de $\mathbf{D}_p(V)$) et on utilise les notations analogues $I_H(\Delta)$.

Nous supposons dans la suite pour que notre problème ait un intérêt que $0 \in I_H(\mathbf{D}_p(V))$, c'est-à-dire que $u$ et $u^*$ sont $\geq 1$.

Soit $t_N(\mathbf{D}_p(V)) = \sum_{r \in \mathbb{Q}} r \dim_K \mathbf{D}_p(V)_r$ où $\mathbf{D}_p(V)_r$ est le sous-espace de $\mathbf{D}_p(V)$ stable par $\varphi$ de pente $r$ (cf. [8]), c'est le degré de Newton du $\varphi$-module filtré $\mathbf{D}_p(V)$.

Le $\varphi$-module filtré $\mathbf{D}_p(V)$ est faiblement admissible au sens de Fontaine, ce qui signifie que

1. $t_H(\mathbf{D}_p(V)) = t_N(\mathbf{D}_p(V))$ ;
2. pour tout sous-$\varphi$-module $\Delta$ de $\mathbf{D}_p(V)$, $t_H(\Delta) \leq t_N(\Delta)$.

De plus, si $\Delta$ est un sous-espace de $\mathbf{D}_p(V)$ stable par $\varphi$ tel que $t_H(\Delta) = t_N(\Delta)$, $\Delta$ est admissible, c'est-à-dire qu'il existe une sous-représentation $V_1$ de $V$ telle que $\mathbf{D}_p(V_1) = \Delta$ ([4, 4.5]).

On dit que $V$ vérifie la propriété ($N_j$) pour un entier $j$ si $V$ n'admet pas de sous-représentation $W$ non nulle telle que $\mathrm{Fil}^j \mathbf{D}_p(W) = 0$. Autrement dit, $V$ vérifie ($N_j$) si et seulement si $V(j)$ vérifie (N). En termes de $\varphi$-modules filtrés, $V$ vérifie ($N_j$) si et seulement si pour tout sous-$\varphi$-module filtré non nul $\Delta$ de $\mathbf{D}_p(V)$ vérifiant $\mathrm{Fil}^j \Delta = 0$, on a $t_H(\Delta) < t_N(\Delta)$ (on dira aussi que $\mathbf{D}_p(V)$ vérifie ($N_j$)).

2.4. **Comportement par produit tensoriel.** Nous aurons besoin de résultats de Totaro sur le comportement des degrés des $\varphi$-modules filtrés par passage au produit tensoriel ([17], voir aussi [15]).

Si $\Delta$ est un $\varphi$-module filtré non nul de dimension $d_\Delta$, posons $\lambda(\Delta) = (t_H(\Delta) - t_N(\Delta))/d_\Delta$. Soit $c \in \mathbb{R}$. On dit qu'un $\varphi$-module filtré $\mathcal{D}$ est de pente $\leq c$ (resp.$< c$) si pour tout sous-$\varphi$-module filtré $\Delta$ non nul de $\mathcal{D}$, on a $\lambda(\Delta) \leq c$ (resp. $\lambda(\Delta) < c$).

Ainsi, un $\varphi$-module filtré faiblement admissible est un $\varphi$-module filtré $\mathcal{D}$ de pente $\leq 0$ tel que $\lambda(\mathcal{D}) = 0$.

**2.4.1. Proposition.** *(Totaro) Soient $\mathcal{D}_i$ pour $i = 1, 2$ deux $\varphi$-modules filtrés de pente $\leq c_i$. Alors $\mathcal{D}_1 \otimes \mathcal{D}_2$ est un $\varphi$-module filtré de pente $\leq c_1 + c_2$.*

**2.4.2. Corollaire.** 1) *Soit $\mathcal{D}$ un $\varphi$-module filtré de pente $\leq c$, alors $\mathcal{D}^{\otimes n}$ et $\wedge^n \mathcal{D}$ sont des $\varphi$-modules filtrés de pente $\leq nc$.*

2) *Soit $\mathcal{D}$ un $\varphi$-module filtré de pente $< 0$, alors $\wedge^n \mathcal{D}$ est un $\varphi$-module filtré de pente $< 0$.*

*Démonstration de* 2). On applique le 1) au sous-$\varphi$-module filtré $\wedge^n \mathcal{D}$ de $D^{\otimes n}$ en prenant pour $c$ le maximum des $\lambda(\Delta)$ pour $\Delta$ sous-$\varphi$-module filtré non nul de $\mathcal{D}$. Comme il y en a un nombre fini, $c$ est strictement négatif. □

On peut réécrire le corollaire de la manière suivante : si $\mathcal{D}$ est un $\varphi$-module filtré tel que $t_H(\Delta) < t_N(\Delta)$ pour tout sous-$\varphi$-module filtré $\Delta$ non nul de $\mathcal{D}$, alors pour tout entier $v$ et tout sous-$\varphi$-module filtré non nul $\Delta'$ de $\wedge^v \mathcal{D}$, on a $t_H(\Delta') < t_N(\Delta')$.

2.5. **Mise en route de la démonstration.** On pose $I_H(\mathbf{D}_p(V)) = \{-u, \ldots u^* - 1\}$. Soit $J$ un sous-ensemble de $I_H(\mathbf{D}_p(V))$. On note $J^c$ le complémentaire de $J$ dans $I_H(\mathbf{D}_p(V))$. Soit $x \in \mathcal{H}(G_\infty) \otimes Z^1_\infty(K, T)$.

**2.5.1. Définition.** Soit $* \in \{e, f, g\}$. On dit ici que $x$ est $(J, *)$-convenable si la projection de $x$ dans $H^1(K_n, V(k))$ appartient à $H^1_*(K_n, V(k))$ pour tout entier $n \geq 0$ et pour tout $k \in J$.

Si $x$ est $(J,g)$-convenable, il existe un élément $f \in \mathbb{Z}_p[[G_\infty]]$, non diviseur de zéro, tel que $fx$ soit $(J,e)$-convenable. Aussi, il suffit en fait d'étudier les éléments $(J,e)$-convenables.

Nous devons rapidement rappeler la théorie de [13]. Nous y avons construit un homomorphisme $\Omega_{V,u}$ de $\mathcal{H}(G_\infty)$-modules $\mathcal{H}_K^{\psi=0} \otimes_K \mathbf{D}_p(V) \to \mathcal{H}(G_\infty) \otimes Z_\infty^1(K,T)$. Grâce à la loi explicite de réciprocité démontrée par Colmez ([2], voir aussi [14, §5]), on peut alors démontrer la proposition suivante : on pose

$$\ell_j = \frac{\log \chi(\gamma)^{-j}\gamma}{\log \chi(\gamma)} = -j + \frac{\log \gamma}{\log \chi(\gamma)}$$

où $\gamma$ est n'importe quel élément d'ordre infini de $G_\infty$.

**2.5.2. Théorème.** *Soit $x \in \mathcal{H}(G_\infty) \otimes Z_\infty^1(K,T)$. Il existe un élément $\mathcal{L}(x) \in \mathcal{H}_K^{\psi=0} \otimes_K \mathbf{D}_p(V)$ tel que*

$$\Omega_{V,u}(\mathcal{L}(x)) = \prod_{j \in I_H(\mathbf{D}_p(V))} \ell_{-j}.x$$

*et on a $\mathfrak{o}_\varphi(\mathcal{L}(x)) = \mathfrak{o}(x) + u^* - 1$ et $\mathfrak{D}_\varphi(\mathcal{L}(x)) = \mathfrak{D}(x) + u^* - 1$. Si de plus $x$ est $(J,e)$-convenable pour $J \subset I_H(\mathbf{D}_p(V))$, il existe $\mathcal{L}_J(x)$ de $\mathcal{H}_K^{\psi=0} \otimes_K \mathbf{D}_p(V)$ tel que*

$$\Omega_{V,u}(\mathcal{L}_J(x)) = \prod_{j \in J^c} \ell_{-j}.x$$

*et on a $\mathfrak{o}_\varphi(g) = \mathfrak{o}(x) + u^* - \sharp J - 1$ et $\mathfrak{D}_\varphi(g) = \mathfrak{D}(x) + u^* - \sharp J - 1$.*

On a bien sûr $\mathcal{L}(x) = \prod_{j \in J} \ell_{-j}.\mathcal{L}_J(x)$. Ici, $\mathfrak{o}_\varphi(g)$ (resp. $\mathfrak{D}_\varphi(g)$) est la borne inférieure des réels $r$ tels que la suite $||p^{rn}(1 \otimes \varphi)^{-n}g||_{\rho_n}$ tend vers 0 (resp. est bornée) lorsque $n$ tend vers l'infini.

On dit que $g$ est à $\varphi$-support dans $\Delta$ si $\Delta$ est un sous-$\varphi$-module de $\mathbf{D}_p(V)$ tel que $g$ appartient à $\mathcal{H}_K \otimes_K \Delta$ et que $\Delta$ est le $\varphi$-support de $g$ si $\Delta$ est le plus petit sous-$\varphi$-module de $\mathbf{D}_p(V)$ vérifiant cette propriété.

Donnons les propriétés de $g = \mathcal{L}_J(x)$ (avec éventuellement $J = \emptyset$). Soit $G \in \mathcal{H}_K \otimes_K \mathbf{D}_p(V)$ tel que $(1 - p^r\Phi)(D^r(G)) = D^r(g)$ pour tout entier $r$ avec $\Phi = \varphi \otimes \varphi$ (quitte à multiplier $g$ par un élément de $\mathbb{Z}_p[[x]]$, $G$ existe). Comme pour $j \in J^c$, $\pi_{n,j}(\Omega_{V,u}(g)) = 0$ pour tout entier $n \geq 0$ où $\pi_{n,j}$ est la projection de $\mathcal{H}(G_\infty) \otimes Z_\infty^1(K,T)$ dans $H^1(K_n, V(j))$, on en déduit par définition de $\Omega_{V,u}$ que $\varphi^{-n}(D^{-j}(G)(\zeta - 1)) \in K_n \otimes_K \mathrm{Fil}^j \mathbf{D}_p(V)$ pour toute $\zeta$ racine de l'unité d'ordre $p^n$ avec $n \geq 1$ ($\varphi$ est l'extension naturelle de $\varphi$ sur $\mathbb{Q}_p(\mu_{p^n}) \otimes_{\mathbb{Q}_p} \mathbf{D}_p(V) = K_n \otimes_K \mathbf{D}_p(V)$). De la relation

$$D^{-j}(G)(\zeta - 1) - p^{-j}\varphi(D^{-j}(G)(\zeta^p - 1)) = D^{-j}(g)(\zeta - 1),$$

on en déduit que pour tout entier $n \geq 2$, on a

$$\varphi^{-n}(D^{-j}(f)(\zeta - 1)) \in K_n \otimes_K \mathrm{Fil}^j \mathbf{D}_p(V).$$

Ainsi, si $x$ est $(J,e)$-convenable, pour tout entier $n \geq 2$ et tout $j \in J^c = I_H(\mathbf{D}_p(V)) - J$, les vecteurs $\varphi^{-n}(D^{-j}(g)(\zeta_n - 1))$ de $K_n \otimes_K \mathbf{D}_p(V)$ appartiennent à un sous-espace $W_n^j = \varphi^n(K_n \otimes_K \mathrm{Fil}^j \mathbf{D}_p(V)) = K_n \otimes_K \varphi^n \mathrm{Fil}^j \mathbf{D}_p(V)$. Les $W_n^j$ forment une filtration décroissante de $K_n \otimes_K \mathbf{D}_p(V)$ de poids de Hodge-Tate les $h_j(\mathbf{D}_p(V))$. Si $\Delta$ est un sous-$\varphi$-module filtré de $\mathbf{D}_p(V)$ et si $g \in \mathcal{H}_K^{\psi=0} \otimes_K \Delta$, le résultat précédent reste valable en remplaçant $\mathbf{D}_p(V)$ par $\Delta$.

On est alors amené à définir le module suivant de fonctions analytiques pour $\mathcal{D}$ un $\varphi$-module filtré quelconque :

**2.5.3. Définition.** Soient $v$ un entier tel que $I_H(\mathcal{D}) \subset ]-\infty, \cdots, v]$ et $J$ un sous-ensemble de $I_H(\mathcal{D})$. On pose

$$\mathcal{A}_{v,J}^{(r)}(\mathcal{D})$$
$$= \left\{ g \in \mathcal{H}_K^{\psi=0} \otimes_K \mathcal{D} \text{ t.q.} \left\{ \begin{array}{l} \mathfrak{D}_\varphi(g) \leq v + r, \\ D^{-j}(g)(\zeta_n - 1) \in K_n \otimes_K \varphi^n \operatorname{Fil}^j \mathcal{D} \text{ pour } j \leq v \text{ et } n \geq 1 \\ D^{-j}(g)(\zeta_n - 1) = 0 \text{ pour } j \in J \text{ et } n \geq 1 \end{array} \right. \right\}$$

Lorsque $v \leq 0$, nous noterons $\tilde{\mathcal{A}}_{v,J}^{(r)}(\mathcal{D})$ le sous-ensemble de $\mathcal{H}_K \otimes_K \mathcal{D}$ défini par les mêmes conditions (on ne suppose plus que $\psi(g) = 0$, comme les entiers $\leq v$ sont négatifs, il n'y a pas de problèmes pour définir $D^{-j}(g)$). Les cas particuliers suivants nous intéressent particulièrement (nous supposerons même de plus ensuite que $v = 0$) :

$$\mathcal{A}_v(\mathcal{D}) = \mathcal{A}_{v,\emptyset}^{(0)}(\mathcal{D})$$
$$= \left\{ g \in \mathcal{H}_K^{\psi=0} \otimes_K \mathcal{D} \text{ t.q.} \left\{ \begin{array}{l} \mathfrak{D}_\varphi(g) \leq v \\ D^{-j}(g)(\zeta_n - 1) \in K_n \otimes_K \varphi^n \operatorname{Fil}^j \mathcal{D} \text{ pour } j \leq v \text{ et } n \geq 1 \end{array} \right. \right\}$$

et pour $k \in I_H(\mathcal{D})$,

$$\mathcal{A}_{v,\{k\}}(\mathcal{D}) =$$
$$\left\{ g \in \mathcal{H}_K^{\psi=0} \otimes_K \mathcal{D} \text{ t.q.} \left\{ \begin{array}{l} \mathfrak{D}_\varphi(g) \leq v \\ D^{-j}(g)(\zeta_n - 1) \in K_n \otimes_K \varphi^n \operatorname{Fil}^j \mathcal{D} \text{ pour } j \leq v \text{ et } n \geq 1 \\ D^{-k}(g)(\zeta_n - 1) = 0 \end{array} \right. \right\}$$

Remarquons que $\mathcal{A}_v(\mathcal{D})$ est stable par multiplication par un élément de $K \otimes \mathbb{Z}_p[[x]]$.

Les considérations précédentes donnent une description en termes de fonctions analytiques des éléments de $\mathcal{H}(G_\infty) \otimes Z_\infty^1(K, V)$ et plus particulièrement des éléments $(J, g)$-convenables. Si $J \subset I_H(\mathbf{D}_p(V))$, on note $Z_{\infty, J, g}^1(K, V)$ le $\Lambda$-module des éléments $(J, g)$ convenables de $Z_\infty^1(K, V)$ et plus généralement $(\mathcal{H}(G_\infty) \otimes Z_\infty^1(K, V))_{J,g}^{(r)}$ le $\Lambda$-module des éléments de $\mathcal{H}(G_\infty) \otimes Z_\infty^1(K, V)$ qui sont d'ordre $\leq r$ et $(J, g)$-convenables (si $J = \emptyset$, on le supprime de la notation). Prenons $v = \sup I_H(\mathbf{D}_p(V))$. On a alors des isomorphismes (après tensorisation par l'anneau total des fractions $\operatorname{Frac}(\Lambda)$ de $\Lambda$) induits par $x \mapsto \mathcal{L}(x)$

$$Z^1_\infty(K,V) \leftrightarrow \mathcal{A}_v(\mathbf{D}_p(V))$$
$$Z^1_{\infty,J,g}(K,V) \leftrightarrow \mathcal{A}_{v,J}(\mathbf{D}_p(V))$$
$$(\mathcal{H}(G_\infty) \otimes Z^1_\infty(K,V))^{(r)} \leftrightarrow \mathcal{A}_v^{(r)}(\mathbf{D}_p(V))$$
$$(\mathcal{H}(G_\infty) \otimes Z^1_\infty(K,V))^{(r)}_{J,g} \leftrightarrow \mathcal{A}_{v,J}^{(r)}(\mathbf{D}_p(V))$$

Plus précisément, si par exemple $x \in Z^1_{\infty,J,g}(K,V)$, il existe $\beta \in \Lambda$ tel que $\mathcal{L}(\beta x) \in \mathcal{A}_{v,J}(\mathbf{D}_p(V))$. On peut déduire de ce qui précède la proposition suivante :

**2.5.4. Proposition.** *Soit $V$ une représentation $p$-adique cristalline de $G_K$ et soit $v = \sup I_H(\mathbf{D}_p(V))$. Soit $k \in I_H(\mathbf{D}_p(V))$ tel que $V$ vérifie la condition $(N_k)$. Le théorème 0.3 pour la représentation $V(k)$ est équivalent à la nullité de $\mathcal{A}_{v,\{k\}}(\mathbf{D}_p(V))$.*

Par twist, on peut d'autre part se ramener au cas où $v = 0$, ce qui nous ferons désormais.

## 3. Cas d'une représentation de dimension 2

Le cas de dimension 2 contient l'essentiel des arguments, aussi est-il intéressant de le faire séparément et indépendamment. On fixe une représentation $p$-adique $V$ cristalline irréductible de dimension 2 dont les poids de Hodge-Tate sont 0 et $u > 0$. On cherche à montrer le théorème 0.3 pour la représentation $V(k)$ avec $-u < k \leq 0$. On pose $\mathcal{D} = \mathbf{D}_p(V)$. Par la proposition 2.5.4, il s'agit de montrer que $\mathcal{A}_{0,\{k\}}(\mathcal{D}) = 0$ (l'hypothèse $k > -u$ ne nous servira pas).

3.1. Montrons que si $\Delta$ est un sous-$\varphi$-module strict de $\mathcal{D}$, $\mathcal{A}_0(\Delta) = 0$. Soit $g \in \mathcal{A}_0(\Delta)$ non nul. Alors, comme $\mathfrak{D}_\varphi(g) \leq 0$, $g$ est d'ordre $\leq r$ avec $-r$ le poids de Newton de $\Delta$. D'autre part, nécessairement $\Delta$ est de poids de Hodge $-u$, ce qui implique que pour tout entier $n$, $D^i(g)(\zeta_n - 1) = 0$ pour toute racine de l'unité d'ordre $p^n$ et pour $0 \leq i \leq u-1$. Donc, $g$ est d'ordre $\geq u$ (en fait divisible par $\log^u(1+x)$), ce qui implique sa nullité car $r < u$ (si $r = u$, la représentation $V$ n'est pas irréductible).

Ainsi, si $g \in \mathcal{A}_0(\mathcal{D})$ est non nul, son $\varphi$-support est nécessairement $\mathcal{D}$.

3.2. Montrons maintenant que $\mathcal{A}_{0,\{0\}}(\mathcal{D}) = 0$. Soit $g \in \mathcal{A}_{0,\{0\}}(\mathcal{D})$ non nul. Soit $F$ le déterminant de $g$ et de $\Phi(g)$ dans une base de $\mathcal{D}$. Quitte à passer au complété $P$ de l'extension maximale non ramifiée de $K$, on peut choisir une base $e_1, e_2$ de $P \otimes_K \mathcal{D} = \mathcal{D}_P$ telle que $\varphi e_i = p^{r_i} e_i$ avec $r_1, r_2$ les poids de Newton de $\mathcal{D}$ et écrire $g = \sum_i g_i e_i$. La condition que $\mathfrak{D}_\varphi(g) \leq 0$ implique que $g_i$ est d'ordre $\leq -r_i$ et donc que $F$ est d'ordre $\leq -r_1 - r_2 = u$. D'autre part, $D^i(F)(\zeta_n - 1)$ est nul pour toute racine de l'unité d'ordre $p^n$ pour tout entier $n$ et $0 \leq j \leq u - 1$. Comme $\log(1+x)$ divise $g$, on vérifie facilement que $D^u(F)(\zeta_n - 1)$ est nul lui aussi pour toute racine de l'unité d'ordre $p^n$ pour tout entier $n$. On en déduit que $F$ est divisible par $\log^{u+1}(1+x)$ et qu'il est donc d'ordre $\geq u + 1$. Donc $F = 0$.

On peut alors écrire que $\Phi(g) = ag$ avec $a$ appartenant au corps des fractions de $\mathcal{H}_K$. L'élément $a$ vaut $p^{r_i}\varphi(g_i)/g_i$ avec $i = 1$ ou $2$. Soit $\mathcal{V} = \det(g, D(g))$ calculé dans la base $e_1, e_2$. En utilisant de nouveau les propriétés de $g$ relatives à la filtration et le fait que $\log(1+x)$ divise $g$, on vérifie que $\mathcal{V}$ est divisible par $\log^u(1+x)$, donc d'ordre $\geq u$. Etudions le comportement de $\mathcal{V}$ sous $\varphi$. On a

$$\varphi\mathcal{V} = p^{-r_1-r_2}\det(\Phi(g), \Phi(D(g)))$$
$$= p^{-1}\frac{\varphi(g_1 g_2)}{g_1 g_2}\det(g, D(g)),$$

d'où

$$\frac{\varphi(\mathcal{V})}{\mathcal{V}} = p^{-1}\frac{\varphi(g_1 g_2)}{g_1 g_2}$$

En utilisant le lemme 2.1.1, on en déduit que

$$\mathfrak{O}(\mathcal{V}) = -1 + \mathfrak{O}(g_1 g_2) \leq u - 1.$$

Donc $\mathcal{V}$ est nécessairement nul, puisqu'il est d'ordre supérieur à $u$.

La théorie du wronskien permet d'en déduire que les composantes de $g$ dans une base de $\mathcal{D}$ sont liées sur le corps fixe par $D$, c'est-à-dire sur $K$ : $\mu_1 g_1 + \mu_2 g_2 = 0$ avec les $\mu_i \in K$ non tous deux nuls, par exemple $\mu_2$. On peut alors écrire $g = \tilde{g}_1 v$ avec $v \in \mathcal{D}$. La relation $\Phi(g) = ag$ implique que la droite $Kv$ est stable par $\varphi$ et donc que le $\varphi$-support de $g$ est en fait une droite. Ce qui n'est pas possible et implique finalement que $g = 0$.

Nous venons de montrer le théorème 0.3 pour $V$ irréductible de dimension 2 telle que $\mathrm{Fil}^1 \mathbf{D}_p(V) = 0$ et $\mathrm{Fil}^0 \mathbf{D}_p(V) \neq 0$.

3.3. Dans le cas où $k < 0$, l'idée de la démonstration est de remarquer que $\sum_j a_j \Phi^j(g)$ est un élément de $\tilde{\mathcal{A}}_0(\mathcal{D})$ pour tout $a_j \in \mathbb{Z}_p[[x]]$. On cherche alors une telle combinaison qui soit de plus divisible par $\log(1+x)$, c'est-à-dire qui appartienne à $\tilde{\mathcal{A}}_{0,\{0\}}(\mathcal{D})$ et on utilise le fait que ce dernier espace est nul (dans ce qui précède, nous n'avons pas utilisé le fait que $\psi(g) = 0$). Pour cela, remarquons que, comme $g \in \mathcal{A}_{0,\{k\}}(\mathcal{D})$, tous les $D^{-k}(\Phi^n(g))$ sont divisibles par $\log(1+x)$ et qu'il existe un entier $v$ tel que les $\Phi^n(g)$ soient combinaisons linéaires à coefficients dans $K \otimes \mathrm{Frac}(\mathbb{Z}_p[[x]])$ des $\Phi^i(g)$ pour $i < v$, ce qui donne beaucoup de relations entre les $D^j(\Phi^i(g))$ pour $j < -k$ et $i < v$ et permet d'exprimer les $D^j(\Phi^i(g))$ pour $0 < j < -k$ et $i < v$ en termes des $\Phi^i(g)$ pour $i < v$ modulo $\log(1+x)$ et d'obtenir une relation entre les $\Phi^i(g)$ modulo $\log(1+x)$. Par exemple pour $k = -1$, on peut écrire $\Phi^2(g) = a\Phi(g) + bg$ avec $a$ et $b \in K \otimes \mathrm{Frac}(\mathbb{Z}_p[[x]])$ et $h = D(a)\Phi(g) + D(b)g$ convient.

Plus précisément, soit donc $g \in \mathcal{A}_{0,\{k\}}(\mathcal{D})$ pour $k < 0$. On pose $g_i = \Phi^i(g)$ pour $i \geq 0$. Supposons que $g$ et $g_1 = \Phi(g)$ ne sont pas liés : on pose

$$g_i = a_i \Phi(g) + b_i g$$

Les $a_i$ et $b_i$ appartiennent à $K \otimes \mathrm{Frac}(\mathbb{Z}_p[[x]])$ : il suffit de le faire pour $a_2$ et $b_2$ et le même raisonnement que dans le paragraphe 3.2 montre que $g \wedge \Phi^n(g) = B_n \log^u(1+x)$ avec $B_n \in K \otimes \mathbb{Z}_p[[x]]$ et $B_1$ est par hypothèse non nul ; comme $b_2 = -\frac{\Phi(g \wedge \Phi(g))}{g \wedge \Phi(g)} = -p^* \frac{\varphi(B_1)}{B_1}$ et que $a_2 = \frac{g \wedge \Phi^2(g)}{g \wedge \Phi(g)} = \frac{B_2}{B_1}$ l'assertion est claire.

La propriété fondamentale est donc que tous les $D^{-k}(g_i)$ sont divisibles par $\log(1+x)$. Ainsi, en posant $F \stackrel{\log}{\sim} G$ si $F - G \in \mathrm{Frac}(\mathbb{Z}_p[[x]]) \otimes \mathcal{H}_K \log(1+x)$, on a pour tout entier $i$, $D^{-k}(g_i) \stackrel{\log}{\sim} 0$. On pose $k' = -k$. Si on dérive $k'$ fois la relation donnant $g_i$, on en déduit que pour tout entier $i \geq 2$,

$$\sum_{j=1}^{k'-1} \binom{k'}{j} \left( D^{k'-j}(a_i) D^j(g_1) + D^{k'-j}(b_i) D^j(g) \right) + D^{k'}(a_i) g_1 + D^{k'}(b_i) g \stackrel{\log}{\sim} 0 \ .$$

ce que nous écrirons

$$\sum_{j=1}^{k'-1} a_i^j D^j(g_1) + b_i^j D^j(g) \stackrel{\log}{\sim} -a_i^{k'} g_1 - b_i^{k'} g \ .$$

Ainsi, si $H$ est la matrice colonne à $2(k'-1)$ lignes

$$\left( D(g_1), \cdots, D^{k'-1}(g_1), D(g), \cdots, D^{k'-1}(g) \right)^t ,$$

$C_r$ la matrice colonne à $r$ lignes $\left( (-a_i^{k'} g_1 - b_i^{k'} g)_{i=2,\cdots,r+1} \right)^t$ et si $A_r$ est la matrice $(r, 2(k'-1))$ dont la $i$-ième ligne (pour $i = 2, \cdots, r+1$) est $\left( (a_i^j)_{j=1,\cdots,k'-1}, (b_i^j)_{j=1,\cdots,k'-1} \right)$, on a

$$A_r.H \stackrel{\log}{\sim} C_r \ .$$

(Les $g_i$ sont des vecteurs, mais le sens de ce qui précède est clair). De l'existence de $H_r$, on en déduit que si $B_r$ est la matrice de type $(r, 2k'-1)$ obtenue en juxtaposant $A_r$ et $C_r$, le déterminant $J$ de la matrice carrée $B_{2k'-1}$ est divisible par $\log(1+x)$. Par construction, $J$ est de la forme $Ug_1 + Vg_0 = U\Phi(g) + Vg$ avec $U$ et $V$ appartenant à $K \otimes \mathrm{Frac}(\mathbb{Z}_p[[x]])$.

Supposons d'abord $U$ et $V$ nuls. Posons $r_0 = 2k' - 1$. Les déterminants de $(A_{2k'-1}, (D^{k'}(a_i))_{i=2,\cdots,r_0+1}$ et de $(A_{2k'-1}, (D^{k'}(b_i))_{i=2,\cdots,r_0+1}$ sont nuls. Soit $s$ le plus petit entier tel que les $s$ premières lignes de $A_{r_0}$ ne soient pas indépendantes (il existe et est inférieur ou égal à $k'$. La $s$-ième ligne s'exprime donc en fonction des $s-1$ premières : il existe des fonctions $\lambda_i$ tels que

$$\begin{cases} D^i(b_{s+1}) = & \sum_{j=2}^{s} \lambda_j D^i(b_j) \\ D^i(a_{s+1}) = & \sum_{j=2}^{s} \lambda_j D^i(a_j) \end{cases}$$

pour $i = 1, \cdots, s$. En dérivant, on en déduit que

$$\begin{cases} \sum_{j=2}^{s} D(\lambda_j) D^i(b_j) = & 0 \\ \sum_{j=2}^{s} D(\lambda_j) D^i(a_j) = & 0 \end{cases}$$

pour $i = 1, \cdots, s-1$. Le rang de $A_{s-1}$ étant $s-1$, on en déduit que les $D(\lambda_i)$ sont nuls et donc que les $\lambda_i$ sont des constantes. D'où une relation du type

$$\begin{cases} \sum_{j=2}^s \lambda_j b_j = & -\lambda_1 \\ \sum_{j=2}^s \lambda_j a_j = & -\lambda_0 \end{cases}$$

avec les $\lambda_i$ constants. On en déduit que

$$\sum_{j=0}^{s+1} \lambda_j \Phi^j(g) = 0$$

avec $\lambda_i \in K$. En appliquant un certain nombre de fois l'opérateur $\psi$ et en utilisant le fait que $\psi(g) = 0$, on en déduit que $g = 0$.

Supposons maintenant $U$ et $V$ non tous deux identiquement nuls. Par construction, ce sont des éléments de $\mathrm{Frac}(\mathbb{Z}_p[[x]])$. Il existe un élément $S$ de $\mathbb{Z}_p[[x]]$ tel que $h = S(U\Phi(g) + Vg) \in \mathcal{H}_K \otimes_K \mathbf{D}_p(V)$ soit divisible par $\log(1+x)$. Il appartient en fait à $\tilde{\mathcal{A}}_{0,\{0\}}(\mathbf{D}_p(V))$ et est donc nul grâce aux résultats de 3.2. On en déduit que $g$ et $\Phi(g)$ sont liés.

On peut donc poser $\Phi(g) = ag$. Un fois montré que $a \in K \otimes \mathrm{Frac}(\mathbb{Z}_p[[T]])$, le même raisonnement montre que $g$ est nul. Pour cela, utilisons $\mathcal{V} = \det(g, D(g))$. Si $\mathcal{V}$ est nul, $g = 0$ et il n'y a rien à montrer. Si $\mathcal{V}$ est non nul, on vérifie encore que $\mathcal{V}$ est d'ordre $\leq u-1$ et divisible par $\log^{u-1}(x)$. Donc $\mathcal{V} = b \log^{u-1}(1+x)$ avec $b \in K \otimes \mathbb{Z}_p[[x]]$ non nul. D'où

$$a^2 = p^{-u+1} \frac{\varphi(b)}{b} \in K \otimes \mathrm{Frac}(\mathbb{Z}_p[[x]]) .$$

On conclut en utilisant le lemme suivant :

**3.3.1. Lemme.** *Soit $a \in \mathrm{Frac}(\mathcal{H}_K)$ tel que $a^N \in K \otimes \mathrm{Frac}(\mathbb{Z}_p[[x]])$ pour un entier $N$. Alors $a \in K \otimes \mathrm{Frac}(\mathbb{Z}_p[[x]])$.*

*Démonstration.* Quitte à multiplier $a$ par un élément de $\mathbb{Z}_p[[x]]$, on peut supposer que $a^N \in K \otimes \mathbb{Z}_p[[x]]$. Montrons d'abord que $a \in \mathcal{H}_K$. A priori, $a = \alpha/\beta$ avec $\alpha$ et $\beta \in \mathcal{H}_K$. Si $\rho < 1$, $\alpha$ et $\beta$ sont convergentes sur le disque $|x| \leq \rho$, en particulier, $\beta$ a un nombre fini de zéros. En utilisant le théorème de préparation de Weierstrass sur le disque $|x| \leq \rho$, on en déduit que $a$ s'écrit dans ce disque $\alpha_1/\beta_1$ avec $\beta_1$ sans zéros. Donc $1/\beta_1$ est analytique sur $|x| < \rho$ et il en est de même de $a$. En faisant varier $\rho$, on en conclut que $a$ appartient à $\mathcal{H}_K$. Alors $||a||_\rho \leq (||a^N||_\rho)^{1/N} \leq 1$ est borné indépendamment de $\rho < 1$, ce qui implique que $||a||_1$ existe et donc que $a \in K \otimes \mathbb{Z}_p[[T]]$. $\square$

Cela termine la démonstration du théorème 0.3 dans le cas d'une représentation irréductible cristalline de dimension 2.

## 4. Démonstration générale

On prend pour $V$ une représentation $p$-adique cristalline telle que $I_H(\mathbf{D}_p(V)) = \{-u, \cdots, 0\}$. Il s'agit alors de montrer le théorème pour $V(k)$ avec $-u < k \leq 0$ lorsque $\mathbf{D}_p(V)$ vérifie la condition (N$_k$) : tout sous-$\varphi$-module filtré $\Delta$ de $\mathcal{D} = \mathbf{D}_p(V)$ tel que $\mathrm{Fil}^k \Delta = 0$ vérifie $t_H(\Delta) < t_N(\Delta)$.

4.1. Comme dans le cas de dimension 2, la première étape est de montrer le théorème pour $V$ et $k = 0$.

**4.1.1. Lemme.** *Soient $\Delta$ un $\varphi$-module filtré de pente $< 0$ t tel que $\mathrm{Fil}^1 \Delta = 0$ et $g \in \tilde{\mathcal{A}}_0(\Delta)$ tel que $g$ et $\Phi(g)$ soient liés sur le corps des fractions de $\mathcal{H}_K^{\psi=0}$. Alors $g = 0$.*

*Démonstration.* On peut supposer que $\Delta$ est le $\varphi$-support de $g$. Posons $\Phi(g) = ag$. Si les $g_i$ sont les composantes de $g$ dans une base convenable de vecteurs propres de $\Delta_P$, on a $\varphi(g_i) = a g_i$. Soit $\mathcal{V} = g \wedge \cdots \wedge D^{d_\Delta - 1}(g)$. On a

$$\Phi(\mathcal{V}) = p^{d_\Delta(d_\Delta-1)/2} a^{d_\Delta} \mathcal{V} .$$

et

$$\frac{\varphi(\mathcal{V})}{\mathcal{V}} = p^{-d_\Delta(d_\Delta-1)/2} \prod_i \frac{\varphi(g_i)}{g_i} .$$

Donc,

$$\mathfrak{O}(\mathcal{V}) \leq -t_N(\Delta) - d_\Delta(d_\Delta - 1)/2 .$$

En utilisant le lemme 2.2, on vérifie que $\mathcal{V}$ est divisible par $(\log(1+x))^{-t_H(\Delta) - d_\Delta(d_\Delta-1)/2}$ et est donc d'ordre $\geq -t_H(\Delta) - d_\Delta(d_\Delta - 1)/2$. Comme $t_H(\Delta) < t_N(\Delta)$, on en déduit que $\mathcal{V} = 0$. Les composantes de $g$ dans une base de $\Delta$ sont donc liées sur $K$. Supposons $g$ non nul et écrivons $g = \sum_{i \in I} g_i v_i$ où les $v_i$ forment un système libre de $\Delta$ et où les $g_j$ sont linéairement indépendants sur $K$. Le cardinal de $I$ est strictement inférieur à $d_\Delta$. Soit $\Delta'$ le $K$-espace vectoriel engendré par les $v_i$. Montrons que $\Delta'$ est stable par $\varphi$. Dans le cas contraire, comme

$$\Phi(g) = \sum_{i \in I} \varphi(g_i) \varphi(v_i) = a \sum_{i \in I} g_i v_i \in \Delta'$$

il existerait des relations non triviales sur $K$ entre les $\varphi(g_i)$ et donc entre les $g_i$. On en déduit que le $\varphi$-support de $g$ est strictement contenu dans $\Delta$, ce qui contredit l'hypothèse faite sur $\Delta$. Donc $g$ est nul.

$\square$

**4.1.2. Lemme.** *Soit $\Delta$ un $\varphi$-module filtré tel que $\mathrm{Fil}^1 \Delta = 0$ et de pente $< 0$. Alors, $\tilde{\mathcal{A}}_0(\Delta)$ est nul.*

*Démonstration.* Soit $g \in \tilde{\mathcal{A}}_0(\Delta)$, non nul. On peut supposer que le $\varphi$-support de $g$ est $\Delta$. Soit $F = g \wedge \cdots \wedge \Phi^{d_\Delta - 1}(g)$. C'est un élément de $\tilde{\mathcal{A}}_0(\wedge^{d_\Delta} \Delta)$ qui est d'ordre $\leq -t_N(\Delta) = -t_N(\wedge^{d_\Delta} \Delta)$ et d'ordre $\geq -t_H(\Delta)$, par le lemme 2.2. Comme $t_H(\Delta) < t_N(\Delta)$, $F$ est nul. Soit $v$ le plus petit entier tel que $(g, \Phi(g), \cdots, \Phi^{v-1}(g))$ forment un système libre sur $\mathcal{K}$. Posons $f = g \wedge \Phi(g) \wedge \cdots \wedge \Phi^{v-1}(g)$. Il est non nul et vérifie $\Phi(f) = af$ avec $a$ dans le corps des fractions de $\mathcal{H}_K$. La proposition 2.4.1 implique que le $\varphi$-support $\Delta' \subset \wedge^v \Delta$ est de pente $< 0$. Il vérifie donc les hypothèses du lemme 4.1.1 (avec $\Delta$ remplacé par $\Delta'$). Donc, $f = 0$, ce qui contredit la définition de $f$. $\square$

4.2. Soit maintenant $g$ un élément de $\tilde{\mathcal{A}}_{0,\{0\}}(\mathcal{D})$. Nous allons construire un nouveau $\varphi$-module filtré $\tilde{\mathcal{D}}$ tel que $\tilde{\mathcal{D}}$ soit de pente $< 0$ et tel que $g \in \tilde{\mathcal{A}}_0(\tilde{\mathcal{D}})$. En utilisant le lemme 4.1.2, on en déduira que $g = 0$, ce qui démontrera le théorème 0.3 pour $V$.

Notons $\tilde{\mathcal{D}}$ le $\varphi$-module $\mathcal{D}$ muni de la filtration suivante :

$$\mathrm{Fil}^j \tilde{\mathcal{D}} = \begin{cases} \mathrm{Fil}^j \mathcal{D} & \text{si } j < 0 \\ 0 & \text{si } j \geq 0 \end{cases}$$

On a donc

$$\tilde{h}_j(\mathcal{D}) = \begin{cases} h_j(\mathcal{D}) & \text{si } j < -1 \\ h_{-1}(\mathcal{D}) + h_0(\mathcal{D}) & \text{si } j = -1 \\ h_j(\mathcal{D}) & \text{si } j \geq 0 \end{cases}$$

On a $t_H(\tilde{\mathcal{D}}) = t_H(\mathcal{D}) - h_0(\mathcal{D})$. On en déduit que $t_H(\tilde{\mathcal{D}}) < t_N(\mathcal{D})$ puisque $h_0(\mathcal{D})$ est supposé non nul. De plus, pour tout sous-$\varphi$-module filtré $\Delta$ de $\mathcal{D}$, soit $h_0(\Delta) \neq 0$, soit $t_H(\Delta) < t_N(\Delta)$. Donc dans tous les cas, $t_H(\tilde{\mathcal{D}}) < t_N(\tilde{\mathcal{D}})$ et $\tilde{\mathcal{D}}$ est de pente $< 0$. D'autre part, comme par hypothèse $g(\zeta_n - 1) = 0$ pour tout entier $n$, $g$ appartient à $\tilde{\mathcal{A}}_0(\tilde{\mathcal{D}})$.

**Remarques.** 1) La même démonstration convient pour $V(-k)$ si $k$ est un poids de Hodge-Tate de $V$ : on considère pour $\tilde{\mathcal{D}} = \tilde{\mathcal{D}}_k$ le $\varphi$-module $\mathcal{D}$ muni de la filtration suivante :

$$\mathrm{Fil}^j \tilde{\mathcal{D}} = \begin{cases} \mathrm{Fil}^j \mathcal{D} & \text{si } j < k \\ \mathrm{Fil}^{k+1} \mathcal{D} & \text{si } j = k \\ \mathrm{Fil}^j \mathcal{D} & \text{si } j > k \end{cases}$$

On a donc

$$\tilde{h}_j(\mathcal{D}) = \begin{cases} h_j(\mathcal{D}) & \text{si } j < k - 1 \\ h_{k-1}(\mathcal{D}) + h_k(\mathcal{D}) & \text{si } j = k - 1 \\ 0 & \text{si } j = k \\ h_j(\mathcal{D}) & \text{si } j > k \end{cases}$$

De nouveau, $g \in \mathcal{A}_0(\tilde{\mathcal{D}}_k)$ et $\tilde{\mathcal{D}}_k$ est de pente $< 0$.

2) La même démonstration montre que $\mathcal{A}_0^{(\epsilon)}(\mathcal{D})$ est nul pour $\epsilon < 0$ (on change alors la pente de Newton). Ainsi, un élément $g$ de $\mathcal{A}_0(\mathcal{D})$ est soit nul soit vérifie $\mathfrak{O}(g) = 0$.

4.3. Prenons maintenant $k \leq -1$ quelconque, $V$ vérifiant la condition $(\mathrm{N}_k)$.

**4.3.1. Lemme.** *On suppose que $V$ vérifie $(N_0)$. Alors $\mathcal{A}_{0,\{k\}}(\mathbf{D}_p(V)) = 0$.*

Nous démontrerons ce lemme plus tard. On en déduit le corollaire suivant :

**4.3.2. Corollaire.** *Soit $W$ la plus grande sous-représentation de $V$ telle que $\mathrm{Fil}^0 \mathbf{D}_p(W) = 0$. Alors $\mathcal{A}_{0,\{k\}}(\mathbf{D}_p(W)) = \mathcal{A}_{0,\{k\}}(\mathbf{D}_p(V))$.*

Déduisons-en si $V$ vérifie ($N_k$), alors $\mathcal{A}_{0,\{k\}}(\mathbf{D}_p(V)) = 0$. On raisonne par récurrence : pour toute représentation $W$ vérifiant $\mathrm{Fil}^1 \mathbf{D}_p(W) = 0$, de dimension strictement inférieure à la dimension de $V$ et vérifiant ($N_j$) pour un entier $j \leq 0$, $\mathcal{A}_{0,\{j\}}(\mathbf{D}_p(W)) = 0$.

Si $V$ vérifie ($N_0$), le lemme 4.3.1 permet de conclure. Sinon, soit $W$ la plus grande sous-représentation de $V$ telle que $\mathrm{Fil}^0 \mathbf{D}_p(W) = 0$. Elle est de dimension $< \dim V$ et on a par le lemme 4.3.2 $\mathcal{A}_{0,\{k\}}(\mathbf{D}_p(W)) = \mathcal{A}_{0,\{k\}}(\mathbf{D}_p(V))$. Si $k_W$ est le plus grand entier tel que $\mathrm{Fil}^{k_W} \mathbf{D}_p(W) \neq 0$, on a nécessairement $k \leq k_W$ car $V$ vérifie ($N_k$) et $W$ vérifie encore ($N_k$). Par hypothèse de récurrence, $\mathcal{A}_{0,\{k\}}(\mathbf{D}_p(W)) = 0$, ce qui termine la démonstration.

*Démonstration du lemme 4.3.1.* On suppose que $V$ vérifie ($N_0$). Soit $g$ un élément de $\mathcal{A}_{0,\{k\}}(\mathbf{D}_p(V))$ supposé non nul. Soit $\Delta$ son $\varphi$-support. Si $\Delta$ n'est pas faiblement admissible (par exemple si $\mathrm{Fil}^0 \Delta = 0$), $\Delta$ est de pente strictement négative et on peut appliquer le lemme 4.1.2 pour en déduire que $g$ est nul.

Supposons maintenant que $\Delta$ est faiblement admissible et que $\mathrm{Fil}^0 \Delta \neq 0$. Il vérfie encore la condition ($N_0$). On a donc $\tilde{\mathcal{A}}_{0,\{0\}}(\Delta) = 0$. Soit $v$ le plus petit entier tel que $g, \cdots, \Phi^v(g)$ forment un système lié. On a donc $\Phi^v(g) = \sum_{i=0}^{v-1} a_v(i) \Phi^i(g)$. On pose alors

$$\Phi^n(g) = \sum_{i=0}^{v-1} a_n(i) \Phi^i(g)$$

pour tout entier $n \geq v$.

**4.3.3. Lemme.** *Les $a_n(i)$ appartiennent à $K \otimes \mathrm{Frac}(\mathbb{Z}_p[[x]])$ pour $i = 0, \cdots, v-1$ et $n \geq v$.*

*Démonstration.* Il suffit de le démontrer pour $n = v$. Soit $F_i = g \wedge \cdots \wedge \overline{\Phi^i(g)} \wedge \Phi^v g$ pour $i = 0, \cdots, v$, le terme surligné étant omis. On a $F_i = \pm a_v(i) F_v$ pour $i = 0, \cdots, v-1$. Par définition de $v$, $F_v$ est non nul. Soit $\Delta'$ le $\varphi$-support de $F$ dans $\wedge^v \Delta$. Posons $\mathcal{V}_i = \det_{\Delta'}(F_i, D(F_i), \cdots, D^{d_{\Delta'}-1}(F_i))$. Alors, on vérifie comme en 3.2 que $\mathcal{V}_i = (\log(1+x))^{-t_H(\Delta') - d_{\Delta'}(d_{\Delta'}-1)/2} B_i$ avec $B_i \in K \otimes \mathbb{Z}_p[[x]]$. Si $B_i$ est nul, $F_i$ est nul, ce qui implique que $a_v(i)$ est nul. Si $B_i$ est non nul, en appliquant $\Phi$, on obtient

$$\pm p^* a_v(i)^{d_{\Delta'}} = \varphi(B_i)/B_i$$

donc $a_v(i)^{d_{\Delta'}} \in K \otimes \mathrm{Frac}(\mathbb{Z}_p[[x]])$ et par le lemme 3.3.1, $a_v(i)$ appartient à $K \otimes \mathrm{Frac}(\mathbb{Z}_p[[x]])$, ce qui démontre le lemme. $\square$

Posons pour simplifier $g_i = \Phi^i(g)$ pour $i \geq 0$ et $k' = -k$. Tous les $D^{k'}(g_i)$ sont divisibles par $\log(1+x)$. Si on dérive $k'$ fois la relation donnant $g_i$, on en déduit que pour tout entier $i \geq v$,

$$\sum_{l=0}^{v-1} \sum_{j=1}^{k'-1} \binom{k'}{j} D^{k'-j}(a_i(l)) D^j(g_l) + \sum_{l=0}^{v-1} D^{k'}(a_i(l)) g_l \stackrel{\log}{\sim} 0 \ .$$

ce que nous écrirons
$$\sum_{l=0}^{v-1}\sum_{j=1}^{k'-1} a_i^j(l)D^j(g_l) \stackrel{\log}{\sim} -\sum_{l=0}^{v-1} a_i(l)^{k'} g_l .$$

Ainsi, si $H$ est la matrice colonne à $v(k'-1)$ lignes $\left((D(g_l),\cdots,D^{k'-1}(g_l))_{l=0,\cdots,v-1}\right)^t$, $C_r$ la matrice colonne à $r$ lignes $\left(\sum_{l=0}^{v-1} a_i^{k'}(l))_{i=v,\cdots,r+v-1}\right)^t$ et si $A_r$ est la matrice de type $(v, 2(k'-1))$ dont la $i$-ième ligne (pour $i = v,\cdots,r+v-1$) est $\left(((a_i^j(\ell))_{j=1,..k'-1})_{l=0,..v-1}\right)$, on a

$$A_r.H \stackrel{\log}{\sim} C_r$$

De l'existence de $H_r$, on en déduit que si $B_r$ est la matrice de type $(r, v(k'-1)+1)$ obtenue à partir de $A_r$ et $C_r$, le déterminant $J$ de la matrice carrée $B_{v(k'-1)+1}$ est divisible par $\log(1+x)$. Par construction, $J$ est de la forme $\sum_{\ell=0}^{v-1} U(\ell)g_\ell = \sum_{\ell=0}^{v-1} U(\ell)\Phi^\ell(g)$ avec les $U(l) \in K \otimes \text{Frac}(\mathbb{Z}_p[[x]])$.

Supposons d'abord que les $U(\ell)$ ne sont pas tous identiquement nuls. Par construction, ce sont des éléments de $K \otimes \text{Frac}(\mathbb{Z}_p[[x]])$. On a donc trouvé (quitte à le multiplier par un élément de $\mathbb{Z}_p[[x]]$) un élément de $\mathcal{H}_K \otimes_K \Delta$ divisible par $\log(1+x)$. Il appartient en fait à $\tilde{\mathcal{A}}_{0,\{0\}}(\Delta)$ et est nul. On en déduit que les $\Phi^\ell(g)$ pour $\ell = 0,\cdots,v-1$ sont liés ce qui est contradictoire avec l'hypothèse faite sur $v$. Donc tous les $U(\ell)$ sont nuls. Posons $r_0 = v(k'-1)+1$. Cela signifie que les déterminants de $(A_{r_0},(D^{k'}(a_i))_{i=v,\cdots,r_0+v-1})$ et de $(A_{r_0},(D^{k'}(b_i))_{i=v,\cdots,r_0+v-1})$ sont nuls. Soit $s$ le plus petit entier tel que les $s$ premières lignes de $A_{r_0}$ ne soient pas indépendantes (il existe et est inférieur à $k'$) : il existe des fonctions $\lambda_i$ tels que

$$D^i(a_{s+v-1}(\ell)) = \sum_{j=v}^{s+v-2} \lambda_j D^i(a_j(\ell))$$

pour $i = 1,\cdots,s$ et pour $\ell = 0,\cdots,v-1$. En dérivant, on en déduit que

$$\sum_{j=v}^{s+v-2} D(\lambda_j)D^i(a_j(\ell)) = 0$$

pour $i = 1,\cdots,s-1$ et pour $\ell = 0,\cdots,v-1$. Le rang de $A_{s-1}$ étant $s-1$, les $D(\lambda_j)$ sont nuls et les $\lambda_j$ sont des constantes. D'où une relation du type

$$\sum_{j=v}^{s+v-2} \lambda_j a_j(\ell) = -\lambda_\ell$$

pour $\ell = 0, \cdots, v-1$ avec les $\lambda_i$ constants. On en déduit que
$$\sum_{j=0}^{s+v-2} \lambda_j \Phi^j(g) = 0$$
avec $\lambda_i \in K$. Cela n'est pas possible : en appliquant un certain nombre de fois l'opérateur $\psi$ et en utilisant le fait que $\psi(g) = 0$ par exemple, on en déduit que $g = 0$. D'où le lemme 4.3.1 □